\documentclass[a4paper,fleqn]{cas-dc}

\usepackage[numbers]{natbib}

\usepackage[ruled,vlined,linesnumbered]{algorithm2e}
\usepackage{tikz,pgf}
\usetikzlibrary{decorations.markings, decorations.pathmorphing, positioning}
\tikzstyle{every node}=[circle, draw,
                       inner sep=0pt, minimum width=3pt]

\newcommand{\mat}[1]{\mathbf{#1}}
\renewcommand{\vec}[1]{\mathbf{#1}}
\newcommand{\id}{\mathbf{I}}
\newcommand{\Kappa}{\mathfrak{K}}
\newcommand{\R}{\mathbb{R}}

\begin{document}
\let\WriteBookmarks\relax
\def\floatpagepagefraction{1}
\def\textpagefraction{.001}
\shorttitle{Efficient state estimation for gas pipeline networks via low-rank approximations}
\shortauthors{N Stahl, N Marheineke}

\title [mode = title]{Efficient state estimation for gas pipeline networks\\ via low-rank approximations}                      
\tnotemark[1]

\tnotetext[1]{This work is supported by the German Federal Ministry for Economic Affairs and Energy, project \emph{MathEnergy -- Mathematical Key Technologies for Evolving Energy Grids}.}

\author[1]{Nadine Stahl}
\cormark[1]
\ead{nadine.stahl@uni-trier.de}

\author[1]{Nicole Marheineke}
\ead{marheineke@uni-trier.de}

\address[1]{Chair of Modeling and Numerics, FB IV -- Mathematics, Trier University, D-54286 Trier, Germany}

\cortext[cor1]{Corresponding author}

\begin{abstract}
In this paper we investigate the performance of projection-based low-rank approximations in Kalman filtering. For large-scale gas pipeline networks structure-preserving model order reduction has turned out to be an advantageous way to compute accurate solutions with much less computational effort. For state estimation we propose to combine these low-rank models with Kalman filtering and show the advantages of this procedure to established low-rank Kalman filters in terms of efficiency and quality of the estimate.
\end{abstract}

%


\maketitle

\section{Introduction}
Efficient state estimation plays an important role for mo\-del-predictive control of gas pipeline networks, \cite{Hesam2015,Gugat2018,Reddy2006}. Gas network models are large-scale systems of high complexity. They consist of nonlinear partial differential equations coupled with algebraic constraints, whose simulation is computationally very expensive. To reduce the effort, model hierarchies have been established \cite{GasHierarchy_Brouwer} and, recently, structure-preserving model order reduction has been developed \cite{EggerMOR,sailer:pp:2020}. The reduced (low-rank) models are in particular mass-conserving and energy-dissipating. 

In state estimation Kalman filtering is popular. The Kal\-man Filter, developed in the 1960s \cite{Kalman}, is still one of the most prominent filtering algorithms because of its optimality property for linear systems. Since it is not suitable for large-scale systems due to computational reasons, low-rank Kalman filters have been developed, see \cite{KFoverview} and references within. They reduce the computational effort by simplifying the costly and memory-demanding calculations of the error covariances in the filtering process. For example, the Compressed State Kalman Filter \cite{KitanidisFirst} reduces the covariance matrices via projection techniques and the Spectral Kalman Filter \cite{SpecKF} by help of a Taylor expansion, whereas the Ensemble Kalman Filter \cite{EnsembleKF} considers a sample covariance by collecting an ensemble of possible states with a Monte-Carlo method. 

In this paper we propose to use a structure-preserving reduced (low-rank) model in Kalman filtering. As our low-rank approximation is obtained from projection, our approach has similarities with the Compressed State Kalman Filter being a projection-based low-rank filter. However, we show that our approach is significantly superior to established low-rank Kalman filters in terms of computational speed and memory demands, while the obtained state estimates are of competitive accuracy. We also investigate the combination of a low-rank model with a low-rank filter.

The paper is structured as follows: After introducing a model hierarchy for gas pipeline networks and the filtering in Sec.~\ref{sec:problem}, we present our Kalman filtering approach with a reduced model and discuss its similarities and differences to the Compressed State Kalman Filter in Sec.~\ref{sec:filters}. A numerical performance study for the low-rank approximations -- in terms of quality and efficiency -- is carried out in Sec.~\ref{sec:numerics}, using an academic benchmark example as well as a real gas pipeline network from western Germany. 

\section{Problem setting} \label{sec:problem}
Modeling a gas pipeline network, we proceed from the nonlinear isothermal Euler equations and establish a model hierarchy by help of linearization and structure-preserving reduction.
\subsection{Hierarchical gas network models}
The gas pipeline network is described by a directed graph $\mathcal{G} = (\mathcal{V}, \mathcal{E})$, where the pipes are represented by the edges $\mathcal{E}$ and the junctions and inlets/outlets by the nodes, in particular $ \mathcal{V}=\mathcal{V}_{I} \cup \mathcal{V}_{B} $, $\mathcal{V}_{I} \cap \mathcal{V}_{B} = \emptyset$. The gas dynamics in a pipe $e$ of length $l^e$ is modeled by the nonlinear isothermal Euler equations for pressure $p^e$ and mass flux $q^e$, i.e., for $e \in \mathcal{E}$, $(x,t) \in [0,l^e]\times [0,T]$,
\begin{align}
a^e\partial_t p^e = -\partial_x q^e, \quad b^e \partial_t q^e = -\partial_x p^e - d^e \frac{|q^e|}{p^e} q^e, \label{eq:iso}
\end{align}
with constant pipe parameters $a^e$, $b^e$, $d^e$ and coupled via Kirch\-hoff's conditions in every junction $v \in \mathcal{V}_{I}$ 
\begin{align*} 
\sum_{e \in \delta_v^-} q^e(l^e,t) &= \sum_{e \in \delta_v^+} q^e(0,t),\\
p^e(l^e,t) = p^{v}(t), e \in \delta_v^+, &\quad  p^e(0,t) = p^{v}(t), e \in \delta_v^-,
\end{align*}
where $\delta_v^+$,  $\delta_v^-$ denote the sets of all topologically ingoing and outgoing edges to $v$.
At inlets/outlets $v \in \mathcal{V}_B$ we prescribe the pressure $p^v$ as 
\begin{align*}
p^v(t) = u^v(t) &= u^v_D(t) + u^v_S(t), \\
\mathrm{d} u^v_S &= \kappa^v (\mu^v - u_S^v)\,\mathrm{d} t + \sigma^v\mathrm{d} W^v_t.
\end{align*}
The time-varying input $u^v$ consists of a determin\-is\-tic part $u^v_D$ and a stochastic one $u^v_S$, where
we particularly model $u_S^v$ by an Ornstein-Uhlenbeck process with constant parameters $\kappa^v$, $\mu^v$, $\sigma^v$ and a standard Wiener process $W^v_t$. The stochastic boundary data allows the incorporation of, e.g., market strategies or the concept of Power-to-Gas \cite{clees2021}. We initialize the transient model with the stationary solution associated to the unperturbed input $u^v=u^v_D(0)$, $v\in \mathcal{V}_B$.
 A linear model variant is obtained when considering a simplified linearized friction term in \eqref{eq:iso} with respective constant $d_l^e$,
\begin{align} \label{eq:liniso}
a^e\partial_t p^e = -\partial_x q^e, \quad b^e \partial_t q^e = -\partial_x p^e - d_l^eq^e.
\end{align}
For an overview on further hierarchical models see, e.g., \cite{GasHierarchy_Brouwer}.

\subsection{Structure-preserving reduced models}
The unperturbed network models can be embedded in a port-Hamiltonian framework and allow for structure-preser\-ving and robust spatial Galerkin approximations. The Galer\-kin approximation is applicable to finite element discretization and projection-based model order reduction. Under mild assumptions on the ansatz spaces, mass conservation and energy dissipation are ensured \cite{EggerMOR,sailer:pp:2020}.

As in \cite{EggerMOR}, we use mixed finite elements and the moment matching method proposed there for the network model  \eqref{eq:liniso}, yielding linear time-invariant descriptor systems of the form
\begin{subequations}\label{eq:fulldescsys}
\begin{align}\label{eq:descriptorsys}
\mat{E}\,\dot{\mathrm{x}} &= \mat{A}\,\mathrm{x} + \mat{B}\,(\mathrm{u}_D + \mathrm{u}_S), \\\label{eq:out}
 \vec{y}  &= \mat{C}\,\mathrm{x},\\ 
\mathrm{d} \mathrm{u}_S &= \mat{\Kappa}(\vec{\mu} - \mathrm{u}_S)\, \mathrm{d}t + \mat{\Sigma}\, \mathrm{d}\vec{W}_t. \label{eq:fulloup}
\end{align}
\end{subequations}
The input functions $\mathrm{u}_D$, $\mathrm{u}_S$ account for the boundary conditions with diagonal matrices $\mat{\Kappa}$, $\mat{\Sigma}\in \R^{|\mathcal{V}_B| \times |\mathcal{V}_B|}$ containing the entries $\kappa^v$ and $\sigma^v$, respectively, and $\vec{\mu} = (\mu^v)_{v\in \mathcal{V}_B}$. The output  $\vec{y}(t) \in \R^{R}$ with certain pressure or flux values corresponds to measurement data in the context of filtering, $t\in[0,T]$. In case of the full order model the state $\mathrm{x}(t)\in \R^N$ represents the space-discrete pressure and flux, its size is determined by the number of pipes and finite elements for each pipe. In the model reduction the system structure is kept. We introduce an orthonormal projection matrix $\mat{V} \in \R^{N\times n}$, $n\ll N$ with the property $\mathrm{x} \approx \mat{V}\, \hat{\mathrm{x}}$, i.e., $\mat{V}^T\mat{V} = \id$, $\id$ identity.  The reduced quantities are indicated by $\hat{\,}$, the system matrices particularly read 
\begin{align}\label{eq:reduction}
\hat{\mat{E}} = \mat{V}^T\mat{E}\mat{V}, \,\, \, \hat{\mat{A}} = \mat{V}^T\mat{A}\mat{V}, \,\,\, \hat{\mat{B}} = \mat{V}^T\mat{B}, \,\,\, \hat{\mat{C}} = \mat{C}\mat{V}.
\end{align}
Note that we preserve a block structure by setting up the reduced spaces for pressure ${V}_p$ and flux ${V}_q$ separately. Imposing a compatibility condition on the spaces, $\partial_x {V}_p={V}_q$, ensures stable reduced models.

\subsection{Stochastically forced model for filtering}
The state estimation is based on a filtering model, for which we use the linear spatially discretized network model in the full or a reduced version. Uncertainties are incorporated via a system noise being modeled as (driving) uncorrelated Gaussian process in \eqref{eq:descriptorsys}.

Let $\vec{x} (t)= (\mathrm{x}^T, \mathrm{u}_S^T)^T(t) \in \R^{N+ |\mathcal{V}_B|}$, $t\in[0,T]$, consider equidistant time points $t_k=k\tau$, $\tau=T/K$, $k=0,...,K$, then our filtering model in time-discrete form reads
\begin{align} \label{eq:transitioneq}
\vec{x}_{k+1} &= \mat{\Phi}\, \vec{x}_k + \mat{\Psi}\, \vec{u}_k + \vec{w}_k\\\nonumber
\mat{\Phi} &= \begin{pmatrix}
\mat{A}_\tau^{-1}(\mat{E}+\tau(1-\theta)\mat{A}) & \tau\mat{A}_\tau^{-1}\mat{B} \\ 
& (\id + \tau\mat{\Kappa})^{-1}\end{pmatrix}  \\\nonumber
\mat{\Psi} &= \begin{pmatrix}
\tau \mat{A}_\tau^{-1}\mat{B} & \\ 
& \tau(\id + \tau\mat{\Kappa})^{-1}\Kappa\vec{\mu}\end{pmatrix}.
\end{align}
with $\mat{A}_\tau =\mat{E}-\tau\theta\mat{A}$.
Here, $\vec{x}_{k}$ denotes the approximation $\vec{x}_{k}\approx \vec{x}(t_k)$, and $\vec{u}_k = ((\theta\mathrm{u}_D(t_{k+1})+(1-\theta)\mathrm{u}_D(t_{k}))^T, \vec{1}^T)^T$ the respective input at time $t_k$ with vector of ones $\vec{1}\in \R^{|\mathcal{V}_B|}$. The underlying time-integration is based on a $\theta$-scheme for \eqref{eq:descriptorsys} and an 
Euler(-Maruyama)-scheme for \eqref{eq:fulloup}. Our different treatment is motivated by the observation that the Euler-scheme is sufficient for the computation of the boundary data, whereas the more sophisticated $\theta$-scheme allows for a better capturing of the dynamic behavior of the flow quantities. The state noise $\vec{w}_k \sim \mathcal{N}(\vec{0}, \mat{Q})$ is a normal distributed centered random variable with constant diagonal covariance matrix $\mat{Q}$ for every time $t_k$, in particular
\begin{align*}
\mat{Q}=\tau \begin{pmatrix} \mat{Z}\mat{Z}^T & \\ & \mat{\Sigma}\mat{\Sigma}^T \end{pmatrix}.
\end{align*}
It results from the system noise with amplitude $\mat{Z}\in \R^{N\times N}$ (diagonal matrix) added to \eqref{eq:descriptorsys} and the scaled Wiener process in \eqref{eq:fulloup}.

The output is assumed to be measurable with suitable devices. To account for measurement errors, noise is added to the output equation \eqref{eq:out}, i.e.,
\begin{align} \label{eq:measurementeq}
\vec{y}_k &= \mat{H}\, \vec{x}_k + \vec{v}_k 
\end{align}
with $\vec{y}_k$ output at time $t_k$ and $\mat{H} = (\mat{C}, \mat{0})$. The discrete time noise $\vec{v}_k\sim \mathcal{N}(\vec{0},\mat{R})$, with constant covariance matrix $\mat{R} \in \R^{R\times R}$, is assumed to be uncorrelated in time and to be component-wisely independent of $\vec{w}_k$. In this study the measurement data $\{\vec{y}_k\}$ for the filtering is provided by the outputs of the nonlinear network model \eqref{eq:iso}.

Using the reduced order model for filtering, the system matrices are replaced by their reduced counter parts, yielding $\hat{\mat{\Phi}}$, $\hat{\mat{\Psi}}$ and $\hat{\mat{H}}$ in \eqref{eq:transitioneq} and \eqref{eq:measurementeq}. Since the input process is not affected by the model order reduction, the reduced covariance of the state noise is given by $\hat{\mat{Q}} = \mat{V}_\vec{x}^T\mat{Q}\mat{V}_\vec{x}$ with
\begin{align}\label{eq:proj}
\mat{V}_\vec{x} = \begin{pmatrix} \mat{V}& \\ & \id \end{pmatrix}.
\end{align} 

\section{Low-rank approximations} \label{sec:filters}

The focus of this work is on Kalman filtering algorithms, since our filtering model is linear. The Kalman Filter (KF) is known to be optimal for linear systems with white noise as state and measurement noise in the sense that it is unbiased and minimizes the error variance \cite{Kalman}. However, it is not applicable to large-scale systems for computational reasons. We propose a filtering approach based on projection-based model order reduction involving low-rank models and discuss its relation to the established Compressed State Kalman Filter \cite{KitanidisFirst, Kitanidis} being a projection-based low-rank Kalman filter.

Kalman filtering consists of a prediction and a correction step. In each of these steps a state estimate and its error covariance are computed. At time point $t_k$, the state and the error covariance matrix of the prediction step are denoted by $\vec{x}_{k|k-1}$ and $\mat{P}_{k|k-1}$, where the first index indicates the current time point $t_k$ and the second index stands for the time up to which measurement data is considered, i.e., $t_{k-1}$. In the correction step additional measurement data of the current time point is taken into account, yielding the state estimate $\vec{x}_{k|k}$ and the corrected error covariance matrix $\mat{P}_{k|k}$. The correction ratio between prediction and correction is given by the so-called Kalman gain $\mat{K}_{k}$, see Algorithm~\ref{algo:kalman}.

\subsection{Kalman Filter using reduced order models}

In the Kalman Filter the computational effort comes from the determination of the error covariances (and the Kalman gains) being fully occupied matrices of the system's size. The matrix evaluations become memory-intensive and expensive for large-scale systems. Establishing a model hierarchy, we derive a low-rank approximation from projection-based model order reduction. Using a low-rank model to estimate the state and to compute the error covariances significantly reduces the effort and makes the Kalman Filter applicable. Prolongating the resulting low-rank estimate to the high-dimensional space yields then the desired state estimate for the large-scale system.

The Kalman Filter applied to the reduced order model (RKF) gives an estimate for the reduced state $\hat{\vec{x}}_{k|k}$ and its associated error covariance matrix $\hat{\mat{P}}_{k|k}$ of small size. To get the quantities $\vec{x}_{k|k}$ and ${\mat{P}}_{k|k}$ of the full order (large-scale) system, we exploit the linearity of the expectation value and the approximation property of the reduced state $\vec{x} \approx \mat{V}_\vec{x}\, \hat{\vec{x}}$. This leads to  $\vec{x}_{k|k} \approx \mat{V}_\vec{x}\, \hat{\vec{x}}_{k|k}$ and $\mat{P}_{k|k} \approx \mat{V}_\vec{x}\mat{\hat{P}}_{k|k}\mat{V}_\vec{x}^T$.

\subsection{Compressed State Kalman Filter}

The idea of the Compressed State Kalman Filter (CSKF) to overcome the computational complexity of the (classical) Kalman Filter is to replace the error covariance matrix with a low-rank approximation. Thus, it is assumed that there exists a constant orthonormal projection matrix $\mat{V}_\mat{P}$, $\mat{V}_\mat{P}^T\mat{V}_\mat{P}=\id$, such that  $\mat{P}_{k|l} \approx \mat{V}_\mat{P} \tilde{\mat{P}}_{k|l} \mat{V}_\mat{P}^T$ holds for both the predicted ($l = k-1$) and the corrected error covariance matrix ($l = k$) with $\tilde{\mat{P}}_{k|l}$ being of much smaller size. For the filter matrices, low-rank versions are precomputed offline, i.e., $\tilde{\mat{\Phi}} =  \mat{V}_\mat{P}^T\mat{\Phi} \mat{V}_\mat{P}$, $\tilde{\mat{H}} = \mat{H}\mat{V}_\mat{P}$ and  $\tilde{\mat{Q}} = \mat{V}_\mat{P}^T\mat{Q}\mat{V}_\mat{P}$, that are then used for the determination of the low-rank error covariance $\tilde{\mat{P}}_{k|k}$ and the respective low-rank Kalman gain $\tilde{\mat{K}}_{k}$ in each time step (cf., Algorithm~\ref{algo:kalman}, lines 4, 5 and 7). The calculations for the state estimate are performed with the full-rank filter matrices. They stay the same as in the (classical) Kalman Filter, except for prolongating the Kalman gain to the high-dimensional space (${\mat{K}}_{k}\approx \mat{V}_\mat{P} \,\tilde{\mat{K}}_{k}$).

\begin{algorithm}[tb] \label{algo:kalman}
\SetAlgoLined
\tcp{initialization}
$\vec{x}_{0|0} = \vec{x}_0, \mat{P}_{0|0} = \mat{P}_0$\;
 \For{$k = 0, \dots, K-1$}{
 \tcp{prediction}
  $\vec{x}_{k+1|k} = \mat{\Phi}\, \vec{x}_{k|k} + \mat{\Psi}\, \vec{u}_k$\; \label{aline:kpred}
  $\mat{P}_{k+1|k} = \mat{\Phi}\mat{P}_{k|k}\mat{\Phi}^T + \mat{Q}$\; \label{aline:kpcov}
\BlankLine\BlankLine
\tcp{correction}
  $\mat{K}_{k+1} = \mat{P}_{k+1|k}\mat{H}^T(\mat{H}\mat{P}_{k+1|k}\mat{H}^T + \mat{R})^{-1}$\; \label{aline:kgain}
  $\vec{x}_{k+1|k+1} = \vec{x}_{k+1|k} + \mat{K}_{k+1}(\vec{y}_{k+1} - \mat{H}\vec{x}_{k+1|k})$\;
  $\mat{P}_{k+1|k+1} = (\id - \mat{K}_{k+1}\mat{H})\mat{P}_{k+1|k}$\;
 }
 \caption{Kalman Filter}
\end{algorithm}

\subsection{Relation between projection-based filtering variants}

The Kalman Filter on a reduced order model and the Compressed State Kalman Filter seem to be very similar at first glance, since they both deal with projection-based low-rank approximations. However, there exist severe differences.

CSKF still uses the large-scale system for the computation of the state estimates, while the reduction only affects the error covariance matrices and thus the Kalman gain.  In contrast, RKF uses a low-rank model for the entire filtering process, including the state computation. This is expected to further improve the computational speed, but probably at the cost of accuracy compared to CSKF. For a numerical study on this point see Sec.~\ref{sec:numerics}.

As far as the projection is concerned, there exist multiple techniques in the field of model order reduction on how to construct the projection matrix $\mat{V}_\vec{x}$ for a low-rank model, e.g., with Moment Matching, Balanced Truncation or Proper Orthogonal Decomposition (see, e.g., \cite{MORoverview} for an overview). Usually, the special structure of the underlying large-scale system (differential, algebraic or additional input equations) is taken into account to achieve the best possible approximation quality. In CSKF, the choice of $\mat{V}_\mat{P}$ is not obvious, as it exclusively acts on a statistical quantity (error covariance matrix). In the following investigations we use $\mat{V}_\vec{x} = \mat{V}_\mat{P}$.

Applying the same projection matrix, one might expect the same low-rank approximations for the error covariances in both filtering variants when starting with the same initialization. But this is in general not the case, as a closer look at their computation in the prediction step reveals. We find \begin{align*}
\hat{\mat{P}}_{k+1|k} &= \hat{\mat{\Phi}}\hat{\mat{P}}_{k|k}\hat{\mat{\Phi}}^T + \hat{\mat{Q}} && \text{for RKF,}\\
\tilde{\mat{P}}_{k+1|k}& = \tilde{\mat{\Phi}} \tilde{\mat{P}}_{k|k} \tilde{\mat{\Phi}}^T + \tilde{\mat{Q}} && \text{for CSKF}.
\end{align*}
While the state noise covariances satisfy $\tilde{\mat{Q}} = \hat{\mat{Q}}$, the filter state matrices differ, $\hat{\mat{\Phi}}\neq \tilde{\mat{\Phi}}$, in case of an underlying implicit time-discretiza\-tion. See, e.g., the first diagonal block matrix for our pipe network model and the projection matrix \eqref{eq:proj},
\begin{align*}
\hat{\mat{\Phi}}_{11}&=(\mat{V}^T({\mat{E}}-\tau\theta{\mat{A}})\mat{V})^{-1} \mat{V}^T({\mat{E}}+\tau (1-\theta){\mat{A}})\mat{V}\\
&\neq \mat{V}^T({\mat{E}}-\tau\theta{\mat{A}})^{-1} ({\mat{E}}+\tau (1-\theta){\mat{A}})\mat{V} = \tilde{\mat{\Phi}}_{11}.
\end{align*} 
In this context, the projection might be interpreted as a model reduction before or after time-discretization for RKF and CSKF, respectively. 

\subsection{Ensemble Kalman Filter}

Additionally to the projection-based approaches, we consider the Ensemble Kalman Filter (EnKF) for comparison reasons in Sec.~\ref{sec:numerics}, as it is probably the most used Kalman filter variant for large-scale systems, \cite{EnsembleKF}. This filtering algorithm avoids the costly evaluation of the algebraic expression of the error covariance matrix in the prediction step (matrix-matrix-matrix product of full occupied matrices of system's size, cf.\ Algorithm~\ref{algo:kalman}, line~\ref{aline:kpcov}). Applying a Monte-Carlo method, the error covariance is numerically approximated on basis of an ensemble of samples of the stochastic filter model. The calculation of Kalman gain and corrected error covariance follows the procedure of the classical Kalman Filter. The state estimate is given as ensemble average.

This procedure usually requires only a very small amount of memory, as not the whole covariance matrix, but only a few samples, whose numbers are usually much less compared to the state dimension, have to be stored. Moreover, there exists observation matrix-free implementations. As for computational time, there is not much gain to be expected from parallelizing the ensemble filter, since collecting all samples after each time step involves a massive overhead. However, the Ensemble Kalman Filter can be combined straightforward with a reduced order (low-rank) model (REnKF), as filter model and filter technique are independent of each other.

 \section{Numerical investigations} \label{sec:numerics} 

We investigate the performance of the low-rank filtering techniques in terms of accuracy of the state estimates and computational effort. As gas pipeline networks we consider an academic example of small size for which the classical Kalman Filter is still applicable as well as a large-scale real network from western Germany. As measurement data we take the mass flux values at the boundary nodes obtained from the nonlinear network model which also serves as reference.  For EnKF we always use $100$ samples in accordance to the recommendations in \cite{EnsembleKF}. All computations are performed in MatlabR2017b on an Intel Xeon with 2.2GHz using 12 Cores. 

\subsection{Academic diamond network}
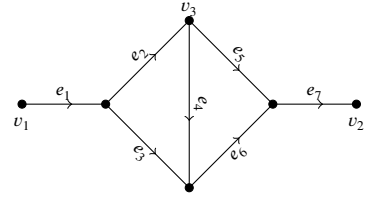
\begin{figure}[tb]
\centering
\tikzset{->-/.style={decoration={ markings, mark=at position 0.6 with {\arrow{>}}},postaction={decorate}}}
\begin{tikzpicture}[scale = 1.1]
	\draw[->-] (0,0) node[fill = black]{} -- (1,0) node[draw=none,fill=none, midway,sloped,above] {\scriptsize $e_1$};
	\draw[->-] (1,0)  -- (2,1) node[draw=none,fill=none, midway,sloped,above] {\scriptsize $e_2$};
	\draw[->-] (1,0) node[fill = black]{}-- (2,-1) node[draw=none,fill=none, midway,sloped,below] {\scriptsize $e_3$};
	\draw[->-] (2,1) -- (2,-1) node[draw=none,fill=none, midway,sloped,above] {\scriptsize $e_4$};
	\draw[->-] (2,1) node[fill = black]{}-- (3,0) node[draw=none,fill=none, midway,sloped,above] {\scriptsize $e_5$};
	\draw[->-] (2,-1) node[fill = black]{} -- (3,0) node[draw=none,fill=none, midway,sloped,below] {\scriptsize $e_6$};
	\draw[->-] (3,0) node[fill = black]{}-- (4,0) node[fill = black]{} node[draw=none,fill=none, midway,sloped,above] {\scriptsize $e_7$};
	\draw[draw=none] (0,-	0.1) node[below, draw=none, fill = none]{\scriptsize $v_1$};
	\draw[draw=none] (4,-	0.1) node[below, draw=none, fill = none]{\scriptsize $v_2$};
	\draw (2,1) node[above , draw=none, fill = none]{\scriptsize $v_3$};
\end{tikzpicture}
\caption{Academic example: diamond network topology.}\label{fig:diamond}
\end{figure}

The small diamond network of 7 edges and 6 nodes serves as academic benchmark setting, see Fig.~\ref{fig:diamond}. The pipe parameters are taken from \cite{clees2021}. For the linear model we use $d_l^e = d^e {|q^e_{av}|}/{p^e_{av}}$, where $p^e_{av}$ and $q^e_{av}$ are the average pressure and mass flux on edge $e$ of the stationary solution associated to the input $u^v=u^v_D(0)$, $v\in \mathcal{V}_B$. In the subsequent example each pipe is discretized with $250$ equidistant finite elements, yielding a full order system with $N=3511$ degrees of freedom. The reduced model is of size $n=29$, which implies a low-rank filter model of size $n+|\mathcal{V}_B|=31$ when adding the two equations for the boundary data. The reduction error is of order $\mathcal{O}(10^{-4})$.
Concerning the filter model, the system noise is exclusively added to the dynamic equations, it is set up with respect to the maximal deviation of the stationary solution from its average values (pressure, mass flux) on each pipe. The algebraic constraints are not perturbed, as they describe mass conservation across junctions. Note that they are preserved under the applied model order reduction. The measurement noise accounts for $1\%$ deviation from the maximal measured flux values.
As input functions we choose
\begin{align*}
u_{D}^{v_1}(t) = \begin{cases} 2+t, & 0 \leq t < 1,\\ 3, & 1 < t < 5,\\ 1.5-0.1t, & 5 \leq t < 10, \\ 2, & t \geq 10, \end{cases} \qquad  u^{v_2}(t) \equiv 2,
\\ (\mu,\kappa,\sigma)^{v_1}=(0,3,0.2)
\end{align*}
for $t\in [0,20]$. We apply $1000$ time steps and set $\theta = 0.51$.

\begin{figure}[tb]
\includegraphics[width=\columnwidth]{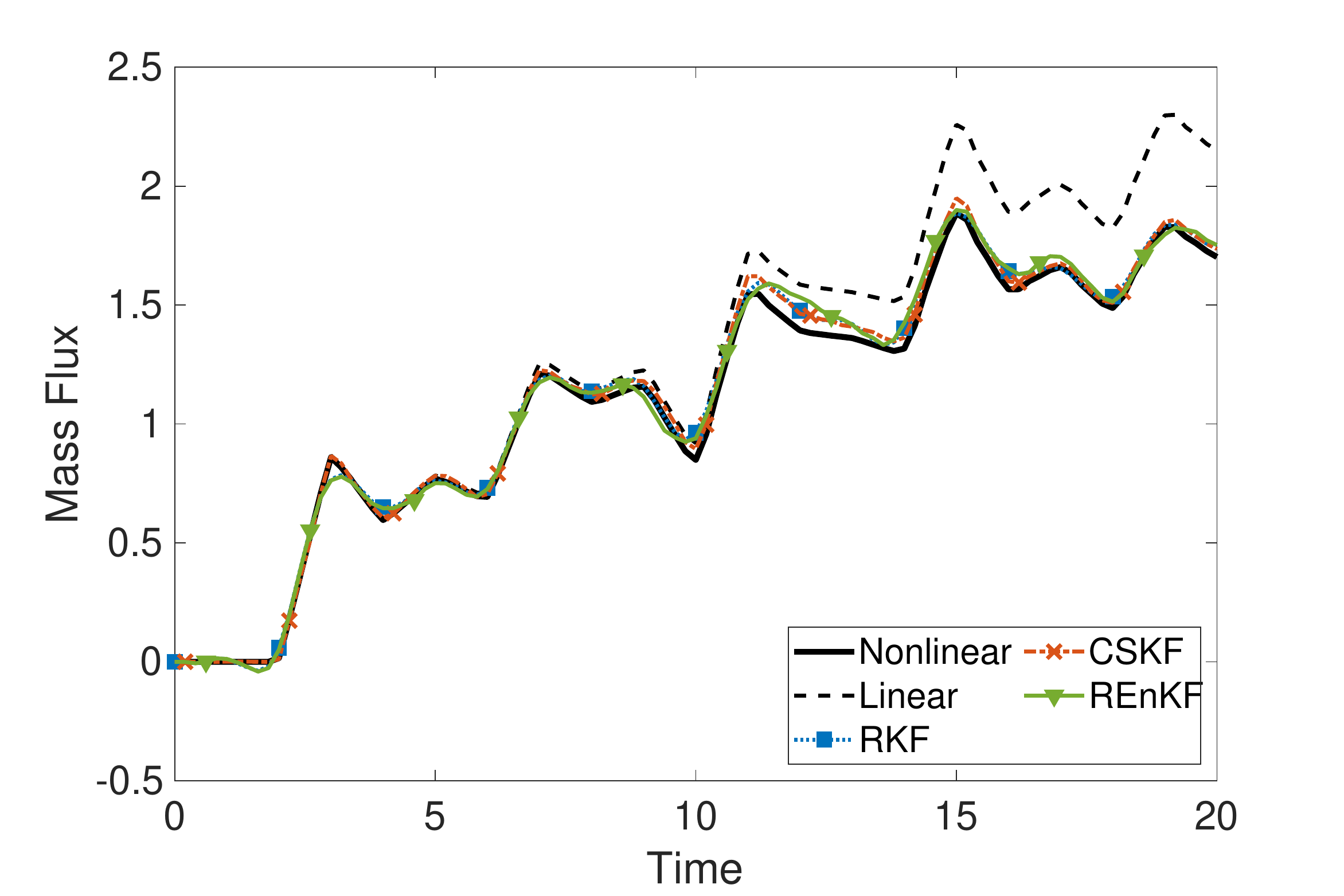}
\caption{Simulation and estimation results for $q^{v_3}$, cf.\ Fig.~\ref{fig:diamond}.} \label{fig:paths}
\end{figure}

\begin{figure}[tb]
\includegraphics[width=\columnwidth]{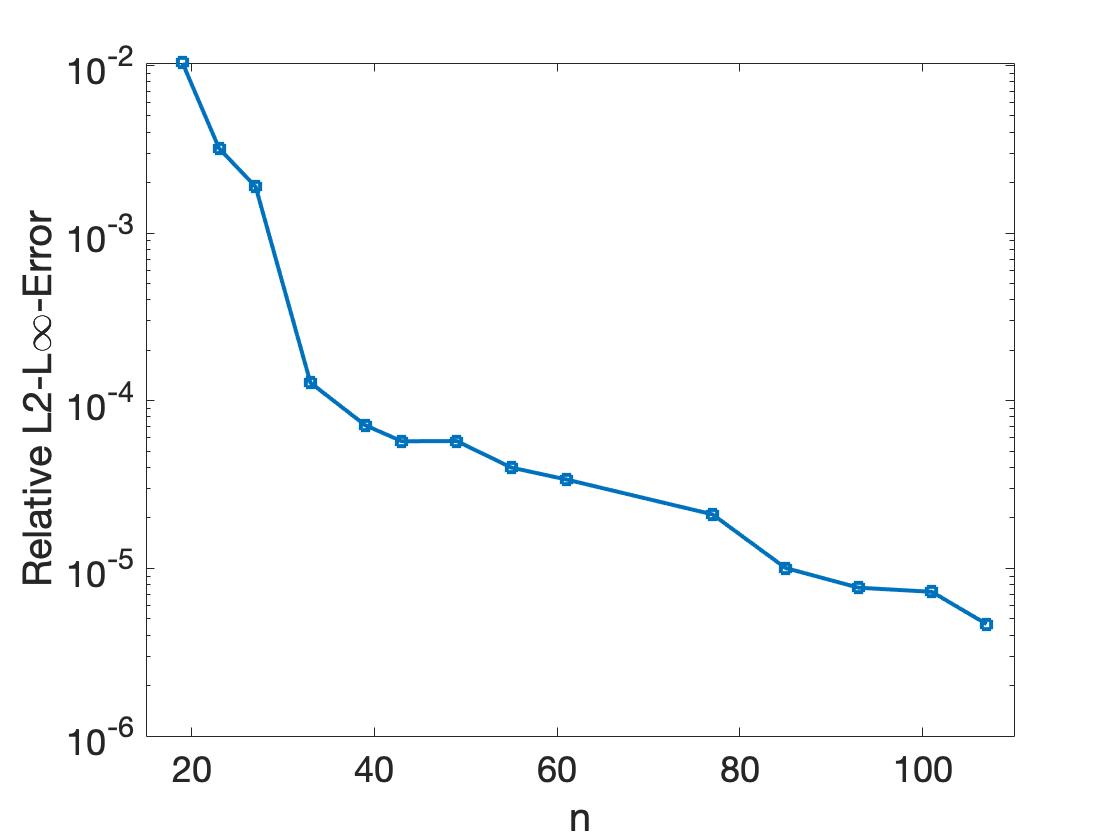}
\caption{Model order reduction error ($L^2$ in space, $L^\infty$ in time) over system size.} \label{fig:morconv}
\end{figure}

Figure~\ref{fig:paths} illustrates the temporal evolution of the mass flux at the inner node $v_3$ (cf.\ Fig.~\ref{fig:diamond}) that is obtained from the nonlinear network model \eqref{eq:iso} (reference solution) and the linear model \eqref{eq:liniso} as well as is estimated by the different low-rank filtering techniques. Although the estimates are computed on basis of the linear model whose solution differs strongly from the nonlinear one, they capture the nonlinear behavior being included via the measurements very well. Regarding the estimation errors (i.e., temporal mean of the relative spatial $L^2$-errors wrt.\ reference) in Table~\ref{tab:erraca}, the (classical) Kalman Filter shows the best approximation properties. CSKF achieves a comparable error, whereas the errors of the other low-rank filters (Kalman and Ensemble Filters on reduced order model) are of same order but slightly higher. The better accuracy of CSKF towards the other low-rank filters might be expected as it estimates the states on basis of the large-scale model. However, all results are astonishingly good regarding the underlying model order reduction error (cf.\ Fig.~\ref{fig:morconv}). 

\begin{table}[tb]
\caption{Errors and CPU time for diamond network}\label{tab:erraca}
\begin{tabular}{l|c|c|c}
Filter & $\text{mean}_j \frac{\|\mathbb{E}[\vec{x}_j - \vec{x}_{j|j}]\|_{L^2}}{\|\mathbb{E}[\vec{x}_j]\|_{L^2}}$ & Offline [s] & Online [s] \\ \hline\hline
KF    & $3.7 \cdot 10^{-2}$ & $8.4 \cdot 10^{2\hphantom{-}}$  & $5.2 \cdot 10^{1\hphantom{-}}$ \\
RKF $^\star$   & $6.0 \cdot 10^{-2}$ & $0.2 \cdot 10^{0\hphantom{-}}$  & $7.4 \cdot 10^{-2}$ \\
CSKF  & $3.8 \cdot 10^{-2}$ & $2.1 \cdot 10^{0\hphantom{-}}$  & $1.5 \cdot 10^{0\hphantom{-}}$  \\
EnKF $^\diamond$ & $8.8 \cdot 10^{-2}$ & $1.2 \cdot 10^{3\hphantom{-}}$                & $1.1 \cdot 10^{3\hphantom{-}}$  \\
REnKF $^\star$$^\diamond$ & $9.8 \cdot 10^{-2}$ & $1.2 \cdot 10^{-1}$                 & $3.2 \cdot 10^{-1}$  \\
\end{tabular}\\
$^\star$ additional post-processing for prolongation is required\\
\hspace*{-0.6cm} $^\diamond$ offline time refers to sampling of random variables
\end{table}

In view of computational time we distinguish between offline (pre-computation) and online (actual runtime) phases.
For the (classical) Kalman Filter as well as its projection-based low-rank variants, the error covariance matrices and the Kalman gains can be precomputed offline. So, the runtime only consists of the state estimation taking into account the actual inputs and measurements. In the ensemble filtering this pre-computation is not possible, since all quantities depend on the chosen samples. Changing input or measurement data hence requires a completely new filtering process. While the ensemble filter is memory-saving, the effort of approximating the sample covariance is approximately the same as running the online phase of KF (or RKF in case of REnKF). It scales here even with the number of samples. Parallelization might be possible, but suffers from the overhead from collecting the sample data in every time step. For sophisticated tuning of the filter to get accurate estimates with fewer samples and reach better running times, we refer to \cite{EnsembleKF, Katzfuss2016}. Note that the offline time listed in Table~\ref{tab:erraca} refers to the sampling of the random variables which is quite expensive in Matlab. Our approach RKF (and also REnKF) clearly outperforms the other filters operating on the large-scale model. However, note that an additional post-processing step, i.e., prolongation of the reduced variables, is required to obtain the estimate for the high-dimensional state and its error covariance.

\subsection{Large-scale real pipeline network}
\begin{figure}[tb]
\includegraphics[trim={3.5cm 4.5cm 4.5cm 3cm}, clip, width=\columnwidth]{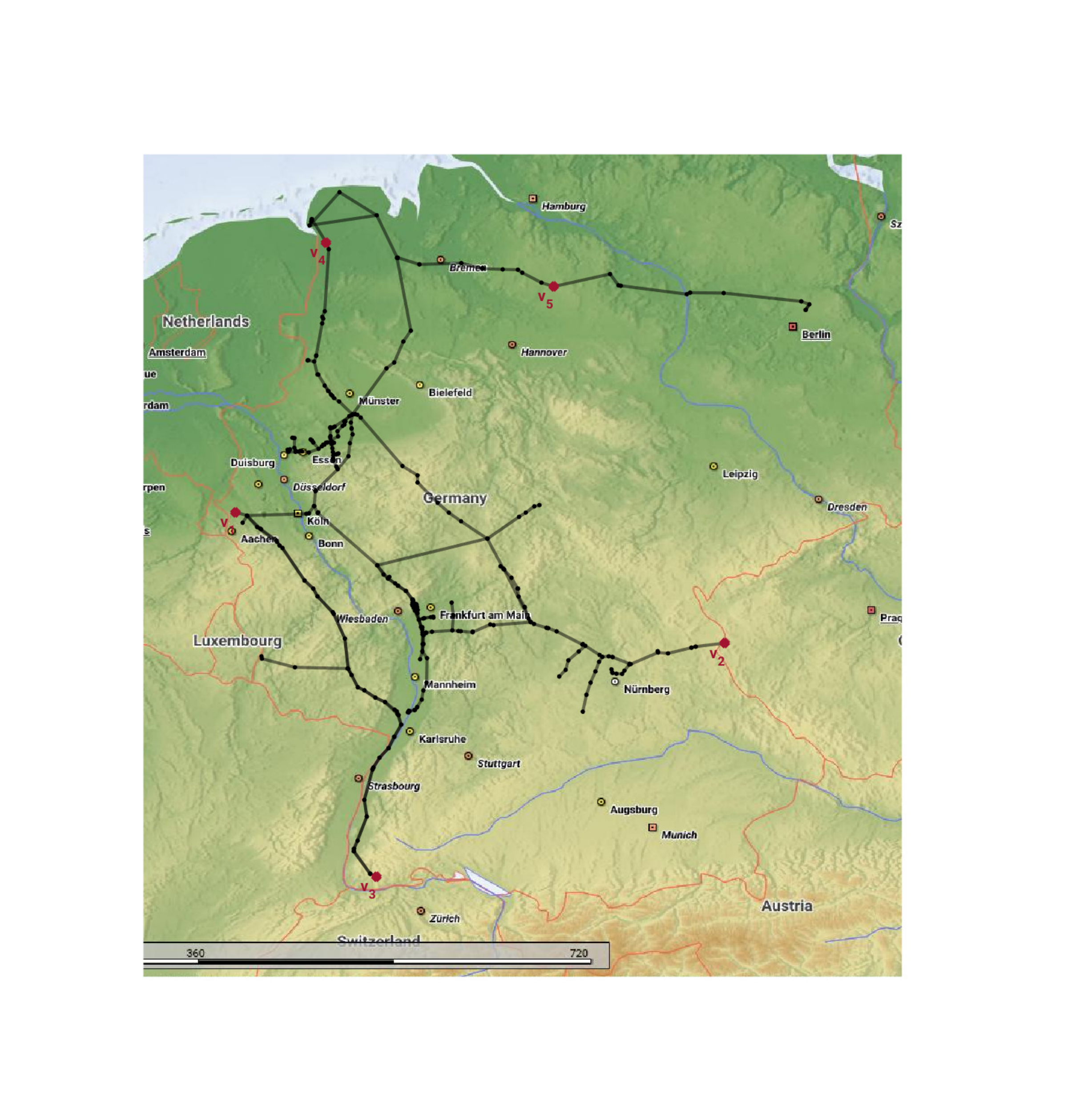}
\caption{Gas pipeline network from western Germany "partDE" \cite{TUBpartDE}.}\label{fig:partDE}
\end{figure}

The partDE network from western Germany consists of $636$ pipes connected in $487$ junctions with $47$ in-/outlets, see Fig.~\ref{fig:partDE}. The pipes have a maximal length of about $120$~km, the specific pipe parameters can be found in \cite{TUBpartDE}. In the subsequent example each pipe is discretized equidistantly with a maximal element length of $100$~m, yielding a full order model of dimension $N = 86160$. The reduced model is of size $n = 1356$, the reduction order is of order $\mathcal{O}(10^{-4})$. As input we prescribe the pressure in bar, particularly we impose 
\begin{align*}
\hspace*{-0.3cm} u_D^{v_1}(t) &= 60 + 5\sin(0.03t), && u_D^{v_2}(t) = 70 + 7\cos(0.1t),\\
\hspace*{-0.3cm} u_D^{v_3}(t) &= 65 + 2\sin(0.05t),\hspace*{-0.3cm} && u_D^{v_4}(t) = 80 + 4\sin(0.008t),\\
\hspace*{-0.3cm} u_D^{v_5}(t) &= 55 + 5\sin(0.017t),\hspace*{-0.3cm}
\end{align*} 
$(\mu,\kappa,\sigma)^{v_i}=(0,3,0.2)$, $i=1,...,5$ at the boundary nodes marked in Fig. \ref{fig:partDE} and $u^{v_j}\equiv60$ at all others, $j=6,...,47$.  A time horizon of $12$~hours is covered by 720 time steps of $1$~minute length. The other model parameters are set and the computations are performed in the same way as in the academic example. 

This real large-scale network makes it necessary to use low-rank approximations for state estimation for computational reasons. The performance of the filters is like in the academic example, see Table~\ref{tab:errun}. Compared to filtering with the reduced order (low-rank) model, the approximation quality of CSKF is slightly better (by less than an order) since the state is estimated on the large-scale model. But, the price to pay for this is significantly longer computing times. An additional drawback, getting more importance the larger the full order model becomes, is the amount of memory needed. This results from the Kalman gains being matrices in dimension of the large-scale model, which are saved due to the pre-computation. The other two methods do not have this drawback, as they operate exclusively on the low dimension and the Kalman gains are hence of small size. The offline phase for the REnKF refers to the sampling or random variables as in the academic example. The time for prolongation of state and error covariance is not included in the listed CPU time as it can be performed for selected time points in a post-processing step.
RKF and REnKF yield comparable results in terms of accuracy and efficiency. For the ensemble filter as a low-rank filtering method operating on a low-rank model, one would expect that the computing costs drop even further. That this is not the case here may be due to the lack of tuning. There might be more sophisticated ways to numerically approximate the covariance matrix from just a few samples, but this is left to future research. However, RKF -- with its separation in offline and online phases -- is clearly preferable, when input or measurement data are changed.

\begin{table}[thpb]
\caption{Errors and CPU time for partDE network}\label{tab:errun}
\begin{tabular}{l|c|c|c}
Filter & $\text{mean}_j \frac{\|\mathbb{E}[\vec{x}_j - \vec{x}_{j|j}]\|_{L^2}}{\|\mathbb{E}[\vec{x}_j]\|_{L^2}}$ & Offline [s] & Online [s] \\ \hline\hline
RKF $^\star$   & $2.2 \cdot 10^{-1}$ & $2.5\cdot 10^{2}$ & $3.4 \cdot 10^{1}$\\
CSKF  & $6.3 \cdot 10^{-2}$ & $2.6\cdot 10^{4}$ & $2.5 \cdot 10^{4}$\\ 
REnKF $^\star$$^\diamond$ & $3.5 \cdot 10^{-1}$ & $7.1 \cdot 10^1$    & $1.4\cdot 10^{2}$
\end{tabular}\\
$^\star$ additional post-processing for prolongation is required\\
\hspace*{-0.6cm} $^\diamond$ offline time refers to sampling of random variables
\end{table}

\section{Conclusion} \label{sec:conclusion}

In this work we presented an efficient approach for state estimation in gas networks. Using a low-rank model derived by linearization and projection-based model order reduction, the Kalman Filter, which is unfeasible for the original large-scale system, becomes applicable. We compared our approach to the Compressed State Kalman Filter whose low-rank also comes from projection-based reduction. Although the underlying approximation idea is very similar, our approach yielding estimates of comparable accuracy proved itself to be much more efficient and less memory-intensive in a numerical study on a German pipeline network. Moreover, the structure-preserving model order reduction can be straightforward combined with other low-rank Kalman filter variants, such as the Ensemble Kalman Filter.

\bibliographystyle{cas-model2-names}

\bibliography{bib}

\end{document}